\documentclass[10pt]{article}

\usepackage{amsmath}
\usepackage{amsthm}
\usepackage{amsfonts}
\usepackage{amssymb}
\usepackage[all]{xy}
\usepackage{tikz}
\usetikzlibrary{arrows}
\usepackage[vmargin=3.5cm]{geometry}

\newcommand{\Z}{\mathbb{Z}}

\def\C{{\mathbb{C}}}

\def\Hom{{\rm{Hom}}}
\def\Ext{{\rm{Ext}}}
\def\Gr{{\rm{Gr}}}

\def\Tr{{\rm{Tr}}}
\def\rep{{\rm{rep}}}

\def\pr{{\rm{pr}}}

\def\dim{{\rm{dim}}\,}
\def\ddim{{\mathbf{dim}\,}}

\def\d{{\mathbf d}}
\def\e{{\mathbf e}}

\def\fl{{\longrightarrow}\,}
\def\ens#1{\left\{ #1 \right\}}
\def\<{\left<}
\def\>{\right>}

\newtheorem{theorem}{Theorem}[section]

\newtheorem{corollary}[theorem]{Corollary}

\theoremstyle{definition}
\newtheorem{remark}[theorem]{Remark}
\newtheorem{example}[theorem]{Example}

\begin{document}

\title{A Homological Interpretation of the Transverse Quiver Grassmannians}
\author{
Giovanni Cerulli Irelli\footnote{Sapienza Universit\`a di Roma, Dipartimento di Matematica. Piazzale Aldo Moro 2, 00185, Roma (ITALY).},
Gr\'egoire Dupont\footnote{Universit\'e de Sherbrooke, D\'epartement de Math\'ematiques. 2500, Boul. de l'Universit\'e, J1K 2R1, Sherbrooke QC (CANADA).}
and
Francesco Esposito\footnote{Universit\`a degli studi di Padova, Dipartimento di Matematica Pura ed Applicata. Via Trieste 63, 35121, Padova (ITALY).}
}

\date{}

\maketitle

\begin{abstract}
	In recent articles, the investigation of \emph{atomic bases} in cluster algebras associated to affine quivers led the second--named author to introduce a variety called \emph{transverse quiver Grassmannian} and the first--named and third--named authors to consider the smooth loci of quiver Grassmannians. In this paper, we prove that, for any affine quiver $Q$, the transverse quiver Grassmannian of an indecomposable representation $M$ is the set of points $N$ in the quiver Grassmannian of $M$ such that $\Ext^1(N,M/N)=0$. As a corollary we prove that the transverse quiver Grassmannian coincides with the smooth locus of the irreducible components of minimal dimension in the quiver Grassmannian.
\end{abstract}

\section{Introduction}\label{Sec:Intro}
	Cluster algebras were introduced by S.~Fomin and A.~Zelevinsky in order to design a combinatorial model for understanding total positivity in algebraic groups and canonical bases in quantum groups \cite{cluster1}. A (coefficient-free) cluster algebra $\mathcal A_Q$ is a certain commutative $\Z$-algebra associated with a quiver $Q=(Q_0,Q_1)$ where $Q_0$ denotes its set of vertices and $Q_1$ its set of arrows. 

	A particular role in the theory of cluster algebras is played by the so-called \emph{atomic bases} (also called \emph{canonically positive bases}). The problem of the existence of atomic bases in arbitrary cluster algebras is still widely open. For cluster algebras associated to affine quivers of type $A_1^{(1)}$ or $A_2^{(1)}$, combinatorial constructions of atomic bases were provided in \cite{shermanz,Cerulli:A21}. 

	In \cite{Dupont:transverse}, it was proved that the elements in these atomic bases can be realised as generating series for Euler-Poincar\'e characteristics of certain constructible subsets of quiver Grassmannians of representations of $Q$, called \emph{transverse quiver Grassmannians} (see Section \ref{ssection:TransGrass} for definitions). However, the definition of transverse quiver Grassmannians is rather technical and can only be formulated for indecomposable representations of affine quivers.

	Independently, the first--named and third--named authors realised the atomic bases in cluster algebras of type $A_1^{(1)}$ as generating series for Euler-Poincar\'e characteristics of the smooth loci in quiver Grassmannians of indecomposable representations of the Kronecker quiver, see \cite{CE:smoothGr}. 

	In this article, we prove that these two realisations are actually the same. More precisely, our main result is the following~:
	\begin{theorem}\label{MainThm:intro}
		Let $Q$ be an affine quiver and let $M$ be an indecomposable representation of $Q$ over $\C$. Then the transverse quiver Grassmannian of $M$ is the set of points $N$ in the quiver Grassmannian of $M$ such that $\Ext^1(N,M/N)=0$.
	\end{theorem}

        The following corollary of Theorem~\ref{MainThm:intro} provides a geometric interpretation of the transverse quiver Grassmannian where $\<-,-\>$ denotes the Euler form on mod-$\C Q$.
	\begin{corollary}\label{CorMainThmIntro}
		Let $Q$ be an affine quiver and let $M$ be an indecomposable representation of $Q$ over $\C$ with dimension vector $\d$. Then for any $\e \in \Z_{\geq 0}^{Q_0}$, the transverse quiver Grassmannian $\Tr_\e(M)$ coincides with the set of smooth points in irreducible components of (minimal possible) dimension $\<\e,\d-\e\>$ in $\Gr_\e(M)$.
	\end{corollary}

	Theorem~\ref{MainThm:intro} and Corollary~\ref{CorMainThmIntro} provide natural interpretations of the transverse quiver Grassmannian whose definition given in \cite{Dupont:transverse} was \emph{ad hoc}. These results also allow one to generalise the definition of the transverse quiver Grassmannians to arbitrary representations of acyclic quivers whereas the initial definition restricted to indecomposable representations of affine quivers. We hope this might be helpful in the construction of atomic bases in a wide class of cluster algebras.

\section{Preliminaries}
	\subsection{Notations and terminology}
		Throughout the article, $Q$ denotes a quiver of affine type, that is an acyclic orientation of an extended Dynkin diagram of type $\widetilde A$, $\widetilde D$ or $\widetilde E$. We denote by $\rep(Q)$ the category of finite-dimensional representations of $Q$ over $\C$, which we identify with the category of finitely generated left-modules over the path algebra $\C Q$ of $Q$. 

		Since $Q$ is representation-infinite acyclic quiver, a representation $M$ of $Q$ can be decomposed uniquely as $M = M_P \oplus M_R \oplus M_I$ where $M_P$ is \emph{preprojective}, $M_R$ is \emph{regular} and $M_I$ is \emph{preinjective}. Moreover, since $Q$ is affine, every regular component of the Auslander-Reiten quiver of mod-$\C Q$ is a \emph{tube} and thus an indecomposable regular representation $R$ of $Q$ can be written uniquely as $R \simeq R_0^{(l)}$ where $R_0$ is \emph{quasi-simple} and $l$ is the \emph{quasi-length} of $R$, see for instance \cite{SS:volume2} for details.

		For any representation $M$ of $Q$, we denote by $\ddim M \in \Z_{\geq 0}^{Q_0}$ its \emph{dimension vector} and given $\e \in \Z_{\geq 0}^{Q_0}$, we denote by
		$$\Gr_{\e}(M) = \ens{N \text{ subrepresentation of } M \ | \ \ddim N = \e}$$
		the \emph{quiver Grassmannian of $\e$--dimensional subrepresentations of $M$}, which is a projective variety. The \emph{quiver Grassmannian} of $M$ is the (finite) union
		$$\Gr(M) = \bigsqcup_{\mathbf e \in \Z_{\geq 0}^{Q_0}} \Gr_{\mathbf e}(M).$$

		A $\C Q$-module $M$ is called \emph{rigid} if $\Ext^1(M,M)=0$. If $M$ is rigid, it is proved in \cite{CR} that for any dimension vector $\e \in \Z_{\geq 0}^{Q_0}$, the quiver Grassmannian is smooth. However, $\Gr_{\e}(M)$ may be non-smooth if $M$ is not rigid.

	\subsection{Transverse quiver Grassmannians}\label{ssection:TransGrass}
		We recall the definition of the transverse quiver Grassmannian provided in \cite{Dupont:transverse}. Let $M$ be an indecomposable $\C Q$-module. If $M$ is rigid, we set $\Tr(M)=\Gr(M)$.

		If $M$ is non-rigid, it is regular and it is thus contained in a tube $\mathcal T$ of rank $p \geq 1$. We denote by $R_i, i \in \Z/p\Z$ the quasi-simple modules in $\mathcal T$ ordered such that $\tau R_i \simeq R_{i-1}$ for any $i \in \Z/p\Z$. Up to relabeling the $R_i$'s, there exist $l \geq 1$ and $0 \leq k \leq p-1$ such that $M \simeq R_0^{(lp+k)}$.

		For any $m \geq n$, $\Hom(R_i^{(n)}, R_i^{(m)})= \C \iota$ for some monomorphism $\iota$ so that we can identify $\Gr(R_i^{(n)})$ with a closed subset of $\Gr(R_i^{(m)})$. By definition, the \emph{transverse quiver Grassmannian} of $M$ is thus~:
		$$\Tr(R_0^{(lp+k)})=\Gr(R_0^{(lp+k)}) \setminus \Gr^{R_0^{(k+1)}}(R_0^{(lp-1)})$$
		where $\Gr^{R_0^{(k+1)}}(R_0^{(lp-1)})$ is the set of subrepresentations in $\Gr(R_0^{(lp-1)})$ containing $R_0^{(k+1)}$ as a subrepresentation. In other words,
		$$\Tr(M)=\Gr(M) \setminus \ens{N \in \Gr(M)\ | \ R_0^{(k+1)} \subset N \subset R_0^{(lp-1)}}.$$
	
		For any indecomposable $\C Q$-module $M$ and any dimension vector $\mathbf e$, we set
		$$\Tr_{\mathbf e}(M) = \Tr(M) \cap \Gr_{\mathbf e}(M).$$
	
\section{Proof of Theorem \ref{MainThm:intro}}
	Fix an affine quiver $Q$ and an indecomposable representation $M$ of $Q$. Then we need to prove that
	\begin{equation}\label{Eq:MainThm}
		\Tr(M)=\{N\in \Gr(M)\ | \ \Ext^1(N,M/N)=0 \}.
	\end{equation}

	If $M$ is rigid, then $\Tr(M)=\Gr(M)$ by definition and it follows from the proof of \cite[Corollary 3]{CR} that $\Ext^1(N,M/N)=0$ for any $N \in \Gr(M)$. Thus \eqref{Eq:MainThm} holds for rigid modules.

	We thus assume that $M$ is not rigid. With notations of Section~\ref{ssection:TransGrass}, we write $M \simeq R_0^{(lp+k)}$ for some $l \geq 1$ and $0 \leq k \leq p-1$ where $p \geq 1$ is the rank of the tube containing $M$.

	In order to prove \eqref{Eq:MainThm} we prove the following two facts:
	\begin{eqnarray}\label{Eq:Inclusion1}
		\textrm{if } N\in \Gr^{R_0^{(k+1)}}(R_0^{(lp-1)}) \textrm{ then } \Ext^1(N,M/N) \neq 0~;\\ \label{Eq:Inclusion2}
		\textrm{if } N\in \Tr(M) \textrm{ then } \Ext^1(N,M/N)=0.
	\end{eqnarray}

	Since $M$ is regular, any submodule $N$ of $M$ is of the form $N=N_R\oplus N_P$ where $N_R$ is regular indecomposable (or zero) and $N_P$ is preprojective. Moreover, the quotient $M/N$ is of the form $M/N=(M/N)_R \oplus (M/N)_I$ where $(M/N)_R$ is regular indecomposable (or zero) and $(M/N)_I$ is preinjective. It follows that
	$$\Ext^1(N,M/N) = \Ext^1(N_P \oplus N_R, (M/N)_R \oplus (M/N)_I) \simeq \Ext^1(N_R, (M/N)_R).$$

	Let us prove \eqref{Eq:Inclusion1}. Let $N\in \Gr^{R_0^{(k+1)}}(R_0^{(lp-1)})$. Since $\Hom(R_0^{(k+1)}, N_P)=0$, it follows that $R_0^{(k+1)}$ is a submodule of $N_R$ and thus
	$$
	N_R=R_0^{(t)}\textrm{ for some }k+1\leq t\leq lp-1.
	$$
	Thus, in the tube containing $M$, the module $N_R$ is contained between $R_0^{(k+1)}$ and $R_0^{(lp-1)}$ in the ray of the tube passing through $M$.

	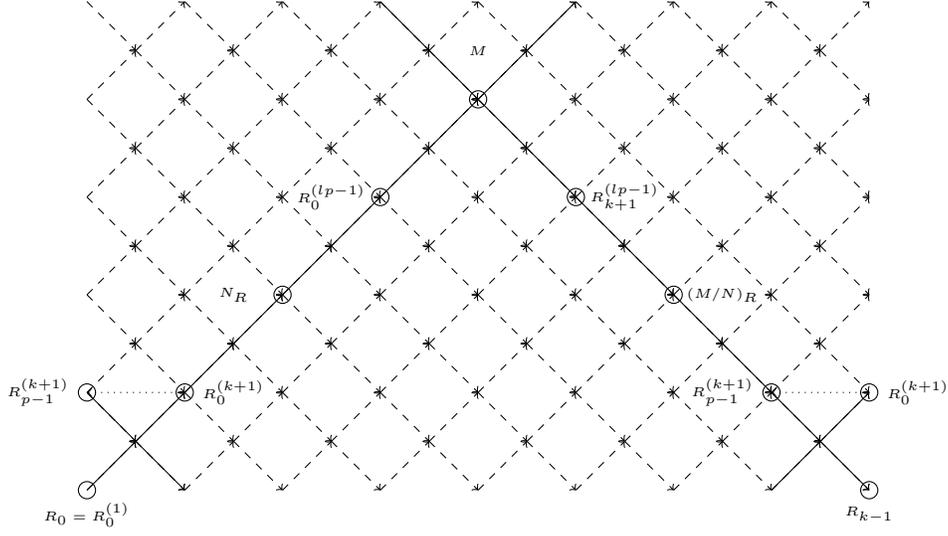
\begin{figure}
		\begin{center}
			\begin{tikzpicture}[scale=.65]
				\tikzstyle{every node}=[font=\tiny]

				\foreach \y in {0,2,4,6,8}
				{
					\foreach \x in {0,2,4,6,8,10,12,14}
					{
						\draw[dashed,->] (\x,2+\y) -- (\x+1,1+\y);
						\draw[dashed,->] (\x+1,1+\y) -- (\x+2,\y);
						\draw[dashed,->] (\x,\y) -- (\x+1,1+\y);
						\draw[dashed,->] (\x+1,\y+1) -- (\x+2,2+\y);
					}
				}
				\draw (0,0) circle (5pt);
				\fill (0,-.5) node { $R_0 = R_0^{(1)}$};

				\draw (0,2) circle (5pt);
				\fill (-1,2) node { $R_{p-1}^{(k+1)}$};

				\draw (2,2) circle (5pt);
				\fill (3,2) node { $R_0^{(k+1)}$};

				\draw (4,4) circle (5pt);
				\fill (3,4) node { $N_R$};

				\draw (6,6) circle (5pt);
				\fill (5,6) node { $R_0^{(lp-1)}$};
				
				\draw (8,8) circle (5pt);
				\fill (8,9) node { $M$};

				\draw (10,6) circle (5pt);
				\fill (11,6) node { $R_{k+1}^{(lp-1)}$};

				\draw (12,4) circle (5pt);
				\fill (13,4) node { $(M/N)_R$};

				\draw (14,2) circle (5pt);
				\fill (13,2) node { $R_{p-1}^{(k+1)}$};

				\draw (16,2) circle (5pt);
				\fill (17,2) node { $R_{0}^{(k+1)}$};


				\draw (16,0) circle (5pt);
				\fill (16,-.5) node { $R_{k-1}$};

				\draw[dotted,->] (16,2) -- (14,2);

				\draw[dotted,->] (2,2) -- (0,2);
	
				\foreach \z in {0,1,2,3,4,5,6,7,8,9}
				{
					\draw[->] (\z,\z) -- (\z+1,\z+1);
					\draw[->] (6+\z,10-\z) -- (6+\z+1,10-\z-1);
				}
				\draw[->] (0,2) -- (1,1);
				\draw[->] (1,1) -- (2,0);
				
				\draw[->] (14,0) -- (15,1);
				\draw[->] (15,1) -- (16,2);

			\end{tikzpicture}
		\caption{The ray and the coray passing through $M=R_0^{(lp+k)}$.}\label{Fig:RayTube}
		\end{center}
	\end{figure}

	Moreover there are surjective morphisms:
	$$
	\xymatrix{
	R_{k+1}^{(lp-1)}\ar@{->>}[r]&(M/N)_R\ar@{->>}[r]&R_{p-1}^{(k+1)}
	}
	$$
	i.e. $(M/N)_R$ is contained between $R_{k+1}^{(lp-1)}$ and $R_{p-1}^{(k+1)}$ in the coray passing through $M$.
	Indeed there are surjective morphisms:
	$$
	\xymatrix{
	M/R_{0}^{(k+1)}\ar@{->>}[r]&(M/N)_R\ar@{->>}[r]&M/R_{0}^{(lp-1)}
	}
	$$
	and clearly
	\begin{equation*}
		M/R_{0}^{(k+1)}\simeq R_{k+1}^{(lp+k-(k+1))}=R_{k+1}^{(lp-1)}
	\end{equation*}
	and
	\begin{equation*}
		M/R_{0}^{(lp-1)}\simeq R_{lp-1}^{(lp+k-(lp-1))}=R_{p-1}^{(k+1)}.
	\end{equation*}
	
	Let $\pi$ be the epimorphism $(M/N)_R \fl R_{p-1}^{(k+1)}$ and let $\iota$ be the monomorphism $\tau R_{0}^{(k+1)} \fl \tau N_R$ (see Figure~\ref{Fig:RayTube}). Then the composition $\iota \circ \pi$ is a non-zero morphism in $\Hom((M/N)_R,\tau N_R)$. It follows from the Auslander-Reiten formula that $\Ext^1(N_R,(M/N)_R) \simeq \Hom((M/N)_R,\tau N_R) \neq 0$. This proves \eqref{Eq:Inclusion1}.

	We now prove \eqref{Eq:Inclusion2}. By definition $N \in \Tr(M)$ if either $R_0^{(k+1)}$ is not a subrepresentation of $N$ or if $N$ is not a subrepresentation of $R_0^{(lp-1)}$.

	If $R_0^{(k+1)}$ is not a subrepresentation of $N$ then $N_R \simeq R_0^{(t)}$ for some $0 \leq t \leq k$. It follows easily that $\Ext^1(N_R,R')=0$ for any regular quotient $R'$ of $R_0^{(lp+k)}$. In particular, we get $\Ext^1(N_R, (M/N)_R)=0$.

	If $N$ is a subrepresentation of $R_0^{(lp+k)}$ but not of $R_0^{(lp-1)}$, we define $a$ to be the minimal integer such that $N \subset R_0^{(lp+a)}$ so that $0 \leq a \leq k$. Consider the composition $f \in \Hom(M,(M/N)_R)$ given by
	$$M \rightarrow M/N \rightarrow (M/N)_R$$
	which is surjective. Since $M$ and $(M/N)_R$ are regular modules and that each tube is closed under kernels, $\ker f$ is a regular module. Since it is a submodule of $M$ and $N \subset \ker f$, we get $\ker f \simeq R_0^{(lp+b)}$ with $a \leq b \leq k$. Thus, there exists an epimorphism
	$$R_0^{(lp+k)}/R_0^{(lp)} \rightarrow R_0^{(lp+k)}/R_0^{(lp+b)} = R_0^{(lp+k)}/\ker(f) \simeq (M/N)_R.$$
	Hence, $(M/N)_R \simeq R_j^{k-j}$ for some $0 \leq j \leq k$ and thus $\Ext^1(R_0^{(s)}, (M/N)_R)$ vanishes for any $s \geq 0$. In particular, $\Ext^1(N_R, (M/N)_R)=0$.
	This finishes the proof of Theorem \ref{MainThm:intro}. \qed

\section{Proof of Corollary \ref{CorMainThmIntro}}\label{section:CorMainThmIntro}
	In this section $Q$ still denotes an affine quiver and $M$ is a representation of $Q$ with dimension vector $\d=(d_i)_{i \in Q_0}$. We denote by $\rep(Q,\d)$ the irreducible affine variety of $\d$-dimensional representations of $Q$. 

	Let $\e=(e_i)_{i \in Q_0} \in \Z_{\geq 0}^{Q_0}$. For any point $N \in \Gr_{\e}(M)$, it is proved in \cite{Schofield:generalrepresentations} that the dimension of the tangent space of $\Gr_{\e}(M)$ at the point $N$ is 
	$$\dim T_N \Gr_{\e}(M) = \dim \Hom(N,M/N).$$
	Thus, it follows from Theorem \ref{MainThm:intro} that $N \in \Tr_{\e}(M) \Leftrightarrow \dim T_N \Gr_{\e}(M) = \<\e,\d - \e\>$.

	For any $i \in Q_0$, we denote by $\Gr(e_i,d_i)$ the set of $e_i$-dimensional sub-vector-spaces of $\C^{d_i}$. Consider the variety
	$$\mathcal G(\e,\d) = \ens{(M,V) \in \rep(Q,\d) \times \prod_{i \in Q_0} \Gr(e_i,d_i) \ | \ V \text{ is a subrepresentation of }M}$$
	and let $\pr_1:\mathcal G(\e,\d) \fl \rep(Q,\d)$ denote the projection onto the first factor. Therefore, for any point $M \in \rep(Q,\d)$, we have $\pr_1^{-1}(M) \simeq \Gr_{\e}(M)$.
	
	Let $X$ denote the closure of the image of $\pr_1$ in $\rep(Q,\d)$. Then $\pr_1$ induces a dominant morphism from $\mathcal G(\e,\d)$ to $X$. Therefore, the dimension of every irreducible component of $\Gr_{\e}(M)$ is greater than $\dim \mathcal G(\e,\d) - \dim X \geq \dim \mathcal G(\e,\d)-\dim \rep(Q,\d) = \<\e,\d-\e\>$, see \cite[II.4, Exercise 3.22 (b)]{Hartshorne}.

	Let $N \in \Gr_{\e}(M)$ and let $Z$ be an irreducible component of $\Gr_{\e}(M)$ containing $N$. Then 
	$$\<\e,\d-\e\> \leq \dim Z \leq \dim T_N \Gr_{\e}(M).$$
	Therefore, $N \in \Tr_{\e}(M)$ if and only if all these inequalities are equalities, that is, if and only if $Z$ has dimension $\<\e,\d-\e\>$ and $N$ is smooth in $Z$. This finishes the proof of Corollary \ref{CorMainThmIntro}. \qed

	\begin{remark}
		Note that we actually proved a stronger result, namely that for any \emph{acyclic} quiver $Q$ and any $\e \in \Z_{\geq 0}^{Q_0}$, the set of points $N \in \Gr_\e(M)$ such that $\Ext^1(N,M/N) =0$ coincides with the set of smooth points in irreducible components of dimension $\<\e,\d-\e\>$ in $\Gr_\e(M)$.
	\end{remark}

\section{Some examples}
	In this section we provide some examples of different phenomena.
	\begin{example}\label{Ex1}
		This example shows that the transverse quiver Grassmannian does not necessarily coincide with the set of smooth points in the quiver Grassmannian in general. Indeed,  let us consider the following acyclic quiver of type $\tilde{A}_{2,1}$
		$$
		\xymatrix@-3ex{
		&&2\ar[dr]&\\
		Q:&1\ar[ur]\ar[rr]&&3,
		}
		$$
		and if $M$ is the indecomposable regular representation
		$$\xymatrix@-3ex{
		&&k^3\ar^{J_3(0)}[dr]&\\
		M=&k^3\ar^{1}[ur]\ar_{1}[rr]&&k^3
		}
		$$
		of dimension $(3,3,3)$ (here $J_3(0)$ denotes the $3\times 3$ indecomposable nilpotent Jordan block). Then, for $\e=(0,2,1)$, it is easy to verify that $\Gr_\mathbf{e}(M)$ is isomorphic to $\mathbb{P}^1$ and is hence smooth of dimension $1 \neq \<\e,\d-\e\>=0$. It is also easily verified that for each submodule $N$ of $M$ with dimension vector $\mathbf e$, we have $\Ext^1(N,M/N) \simeq \C$ so that $\Tr_{\e}(M) = \emptyset$.
	\end{example}

	\begin{example}
		This example points out that we are dealing with quiver Grassmannians viewed as schemes rather than varieties. We hence show an irreducible quiver Grassmannian of dimension $\langle\mathbf{e},\mathbf{d}-\mathbf{e}\rangle$ which does not contain smooth points. Indeed, let $Q$ be the Kronecker quiver, $M$ an indecomposable representation of dimension $\d = (2,2)$ and $\e = (1,1)$. The quiver Grassmannian is the spectrum of the ring of dual numbers. It follows that $\Tr_{\e}(M) = \emptyset$ whereas the variety $\Gr_{\e}(M)$ is a (double) point.
	\end{example}

	\begin{example}
		This example shows a quiver Grassmannian which is not irreducible. Let us consider the same quiver $Q$  of example \ref{Ex1}. Let $M$ be the indecomposable regular representation
		$$\xymatrix@-3ex{
		&&k^2\ar^{J_2(0)}[dr]&\\
		M=&k^2\ar^{1}[ur]\ar_{1}[rr]&&k^2
		}
		$$
		of dimension $(2,2,2)$ and let $\mathbf{e}=(0,1,1)$. The quiver Grassmannian $\Gr_\mathbf{e}(M)$ sits inside $\mathbb{P}^1\times \mathbb{P}^1$ and consists of points $([a:b],[c:d])$ such that $ac=0$. There are two irreducible components, both isomorphic to $\mathbb{P}^1$ which intersects in the point $z=([0:1],[0:1])$. The point $z$ is singular while all the other points are smooth. It is easy to check that the transverse quiver Grassmannian coincides with the set of smooth points.
	\end{example}

\section*{Acknowledgments}
	The first--named author thanks C.~De Concini for the financial support for this project, as part of a post--doctoral scholarship at the Department of Mathematics of ``Sapienza Universit\`a di Roma''. This paper was written while the second--named author was at the university of Sherbrooke as a CRM-ISM postdoctoral fellow under the supervision of Ibrahim Assem, Thomas Br\"ustle and Virginie Charette. The authors would like to thank an anonymous referee for his helpful comments.


\end{document}